[1994/12/01]% LaTeX date must December 1994 or later
\documentclass{article}
\title{Some remarks on the first Hardy-Littlewood conjecture}
\author{Marco Bortolomasi$^1$, Arturo Ortiz-Tapia$^2$}
\date{%
    $^1$Ingegnere e matematico, Ordine degli Ingegneri della Provincia di Modena, c/o Dipartimento di Ingegneria E. Ferrari, Universit\`{a } di Modena e Reggio Emilia, Via P. Vivarelli 10, 41125 Modena; e-mail: bortolamasim@libero.it.\\%
    $^2$PhD in physics Czech Technical University; Universidad Abierta y a Distancia de M\'exico. e-mail: Aortiztapia@nube.unadmexico.mx\\[2ex]%
    \today}

\usepackage{amsmath,amsthm,amsfonts,graphicx}
\usepackage{xcolor}
\usepackage{url}
\usepackage{enumitem}
\usepackage{authblk}

\chardef\bslash=`\\ % p. 424, TeXbook
\theoremstyle{definition}

\theoremstyle{remark}

\newcommand{\eval}[2][\right]{\relax
  \ifx#1\right\relax \left.\fi#2#1\rvert}

\providecommand{\keywords}[1]
{
  \small	
  \textbf{\textit{Keywords---}} #1
}

\begin{document}
\maketitle
\markboth{Some remarks on the first Hardy-Littlewood conjecture}
{Some remarks on the first Hardy-Littlewood conjecture}
\renewcommand{\sectionmark}[1]{}

\abstract{
Starting from the first Hardy-Littlewood conjecture some topics will be covered: an empirical approach to
the distribution of the twin primes in classes mod(10) and a simplified proof of the Brun’s theorem .\\

Finally, it will be explored an approach based on numerical analysis: Monte Carlo Method and Low
discrepancy Sequences will be used to prove the convergence of the conjecture to the expected values.
}


\section{Introduction}

The twin prime conjecture also known as  Polignac’s conjecture is one of the oldest and best-known unsolved problems in number theory and in all of mathematics: it states that for every positive even natural number    $k$, there are infinitely many consecutive prime pairs $p$ and $p'$ such that    $p'-p = k$.\\ The case   $k = 2$ is the twin prime conjecture. Even if the conjecture has not been proved, in spite of many challenges, most mathematicians believe it is true.\\
 Recently, a proof of the conjecture was proposed \cite{arenstorf2004there}, but an error was found after its publication, leaving the conjecture open to this day.\\
What we know for sure, from empirical analysis, is that as numbers get larger, twin primes become increasingly rare.\\
A second twin prime conjecture, called the strong twin prime conjecture or first Hardy-Littlewood conjecture, states that the number $\pi_2(n)$  of twin primes less than or equal to $n$ is asymptotically equal to \footnote{Notation: The use of the asymptotic notations $\mathcal{O}$, $o$, $\sim$ is standard, as well as the symbol $\approx$ used to denote rough, conjectural or heuristic approximations.}:
\begin{equation}\label{EqHLconjecture}
\pi_2(n)\sim 2C_2\int_{2}^{n}\frac{\mathrm{d}x}{(\ln(x))^2}
\end{equation}

where $C_2$ is the so-called twin primes constant \cite{brunconstant2}.\\

Even if both conjectures have not been proved, models for the primes, based on some statistical distribution, can provide the asymptotic value of various statistics about primes. The  ``naive'' Cram\'er random model, models the set of prime numbers by a random set: the starting point is the prime number theorem \footnote{i.e $p(n)\sim n\ln (n)$} involving that in the range $[x,\, x+\varepsilon x]$, for any fixed $\varepsilon >0$ and large $x$, there are about $\frac{\varepsilon x}{\log x}$ primes and each natural number has an independent probability\footnote{but it’s quite obvious that '$p$ is prime'  and '$p+2$ is prime' are not independent events, because $p+2$ is automatically odd and more likely to be prime} of lying in the model set of primes. Using Borel-Cantelli lemma, it can be proved that the model leads to a conjecture of the form:
\begin{equation}
\pi_2(n)\sim\frac{x}{(\ln\ln(x))^2}
\end{equation}
and consequently that there are infinitely many twin primes. The model is too simplified to give accurate results, but tends to give predictions of the right order of magnitude \cite{pintz2007cramer,tao2015} .\\

It is worth noting that in 1996 it was proved \cite{ribenboim2012new} that:
\begin{equation}
\pi_2(n)\leq c \Pi_2\frac{x}{(\ln\ln(x))^2}\left[1+\mathcal{O}\left(\frac{\ln_2(x)}{\ln(x)}\right)\right]
\end{equation}
where $\Pi_2$ is the twin primes constant and $c$ is another constant, that according to Hardy-Littlewood conjecture is 2 and that has been precised to be 6.8325 \cite{haugland1998application} from previous values \cite{twinprimeWeisstein}.

\section{Hardy-Littlewood conjecture: an asymptotic distribution of twin primes}

Let: \\

$$P=\text{ the set of primes}$$

\begin{equation}
X_i(2,m):=\# \{ (p_i,\; p_{i+2}):p_i,\;p_{i+2}\in P \text{ and } \left\lbrace \frac{p_i}{10} \right\rbrace =m \}
\end{equation}
With $i=1,\,2,\,3\; \text{i.e. } X_1(2,1)\;X_2(2,7)\;X_3(2,9)\;$
It is evidently clear that every pair of twin primes, with the sole exception of: $(3\, 5)\,\, (5\, 7)$, belongs to one
of these three classes:
\begin{align*}
X_1(2,1)=& \{(11,\,13)\;(41,\,43)\;(71,\,73)\cdots\}\\ 
X_2(2,7)=& \{(17,\,19)\;(107,\,109)\;(137,\,139)\cdots\}\\ 
X_3(2,9)=& \{(29,\,31)\;(59,\,61)\;(149,\,151)\cdots\}\\ 
\end{align*}
The Hardy-Littlewood conjecture refers to the number of twin primes and doesn't provide any information about their distribution. On the basis of numerical evidence it is possible to propose a different perspective of the famous conjecture, and a correlation, otherwise lacking, between the distribution and the counting function of the twin primes.\\
The distribution of pairs of primes has been studied with the Chi-square $\chi^2$ statistic approach \cite{bortolamasimodello},  in order to compare experimental data to the expected values: based on this analysis, it was possible to verify the hypothesis that twin primes thin out in the three classes with the same cardinality.\\

Let\footnote{$S=\# \bigcup_i\bar{X}_i(2,m)$ differs from $2n$ because classes of only one element (pairs $(3, 5) (5, 7)$) have not been considered in the numerical model}: \\

$\pi_2^i(n)$ be counting function of class $X_i(2,m)$ (e.g. $n=80,\,\pi_2^1(80)=3$).\\
Numerical analysis provides the following results concerning the class $X_1(2,1)$:

%%%%
\begin{figure}[htb]
	\centering
	\includegraphics[width=1\textwidth]{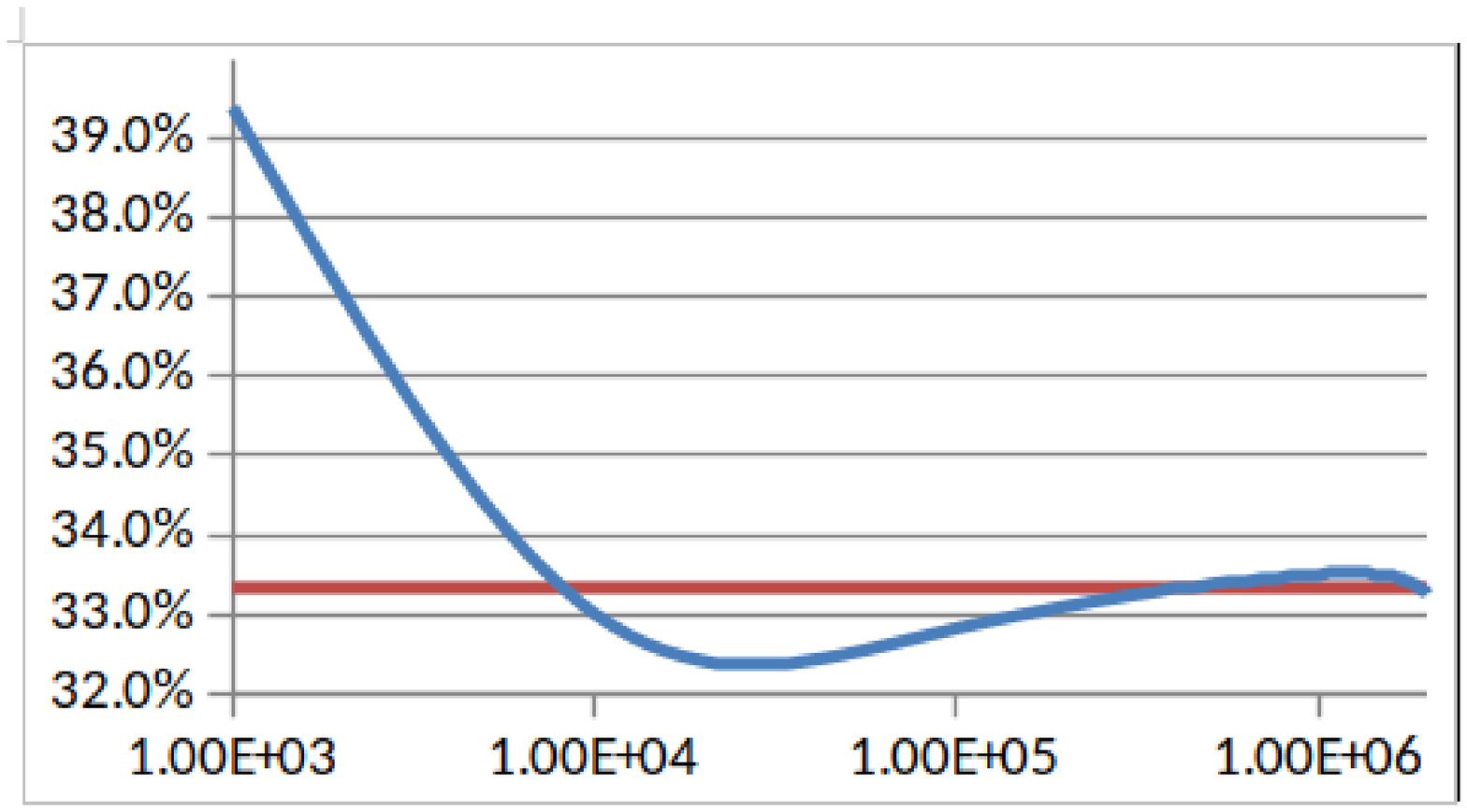}
	\caption{Proportion of $\pi_2^1(n)/\pi_2(n)$ for  $X_1(2,1)$}\label{FigConstants}           
\end{figure}

And similar results for the classes $X_2(2,7)$ and $X_3(2,9)$. It may be clearly seen that the three classes converge toward the same value: 33.3\% and the Chi-square $\chi^2$ statistic approach justifies a random distribution of the twin primes in the three classes.\\

Hence, under empirical evidence, the first Hardy-Littlewood conjecture may be re-written as follows:
\begin{equation}
\pi_2^i(n)=\frac{1}{3}\pi_2(n)\sim \frac{2}{3}C_2\int_{2}^{n}\frac{\mathrm{d}n}{(\ln\ln(n))^2},\;i=1,2,3
\end{equation}

In other words\footnote{It's worth noting that the numerical analysis leads to a similar result also in case of cousin primes, sexy primes and Sophie Germain primes \cite{bortolamasimodello}}, the asymptotic distribution of pairs of twin primes $(p_i,\, p_i + 2)$ in the three classes $X_i(2,m),\,m=p_i\mathrm{mod}(10),\, m=1,\,7,\,9$ may be described as statistically random\footnote{even if some small differences appear in the speed of convergence rate}: no strong empirical evidence appears to the contrary.\\

The fact that twin primes behave more randomly than primes, is also supported by the works by Kelly and Pilling \cite{kelly2001characterization}, \cite{kelly2001discrete} pointing out that the occurrences of twin primes in any sequence of primes are like fixed probability random events.

\section{From Hardy-Littlewood conjecture to the Brun’s theorem}
Viggo Brun wanted to analyze the sum
\begin{equation}\label{EqSumBrunAnalysis}
\sum_{p,\;p+2\;primes} \frac{1}{p}+\frac{1}{p+2}
\end{equation}
  hoping that the sum would be infinite and thus giving a solution to the twin prime conjecture. However, what he proved in 1919, by means of a specific sieve, is that the sum of reciprocals of the twin primes converges to a finite value \cite{brunconstant}.\\
\begin{equation}
\sum_{p,\; p+2 \text{ primes}} \frac{1}{p}+\frac{1}{p+2}\approx 1.9< +\infty
\end{equation}
If the series had diverged, it would have indicated that there is an infinite number of twin primes but the proof that it converges does not provide more information about Polignac’s conjecture. The original proof of the convergence was based on the Brun’s simple pure sieve (principle of Inclusion-Exclusion), although it is possible to provide a simplified demonstration starting from the first Hardy-Littlewood conjecture.

\emph{Proof}\\

First of all, it is easy to observe that:
\begin{equation}
\int_{2}^{n}\frac{\mathrm{d}x}{(\ln\ln (x))^2}\sim \frac{n}{(\ln\ln (n))^2}
\end{equation}

In fact, let:\\
\begin{equation}
f(n)=\int_{2}^{n}\frac{\mathrm{d}x}{(\ln\ln (x))^2} 
\end{equation}
   and 
\begin{equation}
g(n)=\frac{n}{(\ln\ln (x))^2}
\end{equation}
then
\begin{equation}
\frac{f(n)}{g(n)}=\frac{f'(n)}{g'(n)}=\frac{1}{1-2/\ln(n)}=1
\end{equation}

Unfortunately, the asymptotical equivalence does not provide any information about the behavior of the ratio: 
\begin{equation}
\frac{\int_{2}^{n}\frac{\mathrm{d}x}{(\ln\ln (x))^2} }{\frac{n}{(\ln\ln (n))^2}}\in [2,\; +\infty[
\end{equation}

In order to bound the integral with a degree of approximation, in the set $[2,\; +\infty[$  we proceed as follows:
\begin{equation}
\int_{2}^{n}\frac{\mathrm{d}x}{(\ln\ln (x))^2}=\int_{2}^{n}\frac{\mathrm{d}x}{\ln (x)}-\frac{x}{\ln(x)}\Bigg\vert_{2}^{n}=li(n)-li(2)-\frac{n}{\ln(n)}+\frac{2}{\ln(2)}
\end{equation}             
with
\begin{equation}
li(n)=\int_{0}^{n}\frac{\mathrm{d}n}{\ln\ln(x)}
\end{equation}  
The asymptotic expansion (Poincar\'{e} expansion) of $li(n)$ for $n\to \infty$  gives:
\begin{equation}
li(n)\sim \frac{n}{\ln\ln(n)}\sum_{k=0}^{\infty}\frac{k!}{(\ln\ln(n))^k}
\end{equation}
    i.e.\footnote{This implies also:  $li(n)-n\ln(n)=\mathcal{O}(n\ln 2n$)}
    
\begin{equation}
li(n)\sim \frac{n}{\ln(n)}+\frac{n}{\ln^2(n)}+\frac{2n}{\ln^3(n)}+\cdots
\end{equation}        
Hence assuming the Hardy-Littlewood conjecture (Eq.\ref{EqHLconjecture}):
\begin{equation}
\pi_2(n)\sim 2C_2\cdot \left(-li(2)+\frac{2}{\ln(2)}+\frac{n}{\ln^2(n)}+\frac{2n}{\ln^3(n)}+\frac{6n}{\ln^4(n)}+\cdots\right)
\end{equation}

Where $li(2) = 1.045163\cdots$ \cite{li2}
The series is not convergent and an approximation is reasonable where the series is truncated at a finite number of terms with an error roughly of the same size as the next term.\\
  
In fact, the problem associated to divergence is that for a fixed $\varepsilon$, the error in a divergent series will reach to an $\varepsilon$-dependent minimum, but as more terms are added the error then increases without bound and tends to infinity.\\

Since for every $n\in N,\; n \geq 10^{12}$, we have:
\begin{equation}
\frac{1}{\ln^3(n)}\geq \frac{6}{\ln^4(n)}
\end{equation}
   
Hence we can write for every $n\in N,\; n \geq 10^{12}$ i.e. in the set  $[10^{12},\; +\infty[$
\begin{equation}
1< \frac{\pi_2(n)}{2C_2\frac{n}{\ln^2(n)}}\leq 1+\frac{2}{\ln(n)}+\frac{7}{\ln^2(n)}
\end{equation}
i.e.
\begin{equation}
1< \frac{\pi_2(n)}{\frac{n}{\ln^2(n)}}\leq \approx 1.4277
\end{equation}

Hence if we assume the  Hardy-Littlewood conjecture we can say that a number exists $\bar{n}\in N$ such that for every $n\geq \bar{n}$:

\begin{equation}
\pi_2(n)\leq K \frac{n}{\ln^2(n)}
\end{equation}

It is worth noting that the ratio
\begin{equation}
\frac{\pi_2(n)}{\frac{n}{\ln^2(n)}}
\end{equation}

   has been studied by many authors under the general condition:
\begin{equation}
\frac{\pi_2(n)}{\frac{n}{\ln^2(n)}}<2C_2+\varepsilon
\end{equation}

Recentely, Wu \cite{wu2007chen} proved that for a sufficiently large $n$:    

\begin{equation}
\frac{\pi_2(n)}{\frac{n}{\ln^2(n)}}<4.5
\end{equation}

Now let us consider the sum in Eq.\ref{EqSumBrunAnalysis}
\begin{equation}
\sum_{p,\;p+2\;primes} \frac{1}{p}+\frac{1}{p+2}
\end{equation}
Since
\begin{equation}
\frac{1}{p}+\frac{1}{p+2}\leq \frac{2}{p},
\end{equation}
the convergence of Eq.\ref{EqSumBrunAnalysis} is equivalent to the convergence of 
\begin{equation}
\sum_{p,\;p+2\;primes} \frac{1}{p},
\end{equation}
  
there are two possibilities:

\begin{enumerate}[label=\alph*)]
\item Twin primes are finite in number (in this case the sum of the series is finite and the convergence is proved);
\item Twin primes are not finite in number, in this case:
\end{enumerate}

Let $r$ be the $r^{th}$  twin prime\footnote{This part of the proof is the same as in \cite{languasco2007note}} (e.g $q_r=107,\text{ hence } r=\pi_2(107)=10$):
\begin{equation}
r=\pi_2(q_r)\leq K\frac{q_r}{\ln^2(q_r)}\leq K\frac{q_r}{\ln^2(r+1)},\; \text{since } q_r > r+1,\; \forall r\in N
\end{equation}
Hence
\begin{equation}
\frac{1}{q_r}\leq K\frac{1}{r\ln^2(r+1)}
\end{equation}
And:
\begin{equation}
\sum_{p,\;p+2\;primes} \frac{1}{p}=\sum_{1}^{\infty}\frac{1}{q_r}\leq K\sum_{1}^{\infty}\frac{1}{r\ln^2(r+1)}
\end{equation}
For the comparison test, also the series 
\begin{equation}
\sum_{p,\;p+2\;primes} \frac{1}{p}
\end{equation}
converges.

\section{Calculation of the integral $2C_2\int_{2}^{n}\frac{\mathrm{d}x}{(\ln\ln(x))^2}$  using MonteCarlo approach}

Monte Carlo (MC) and Quasi-Monte Carlo (QMC) methods are widely used in numerical analysis, especially in Physics and Finance. Consider an integral of the form: $I=\int_{\Omega} f(x)\mathrm{d}x$. Where $\Omega$ is the domain of integration and $f(x)$ a bounded real function.\\

Most direct quadrature methods are based on the Riemann definition of an integral (a finite sum of ordered '\emph{areas}' under the curve $y=f(x)$): MC and QMC methods are explained by Lebesgue integration: the finite sum do not depend on the order, it is enough that the function can be somehow '\emph{measured}'.\\

By the strong law of large numbers, if $U$ is a uniformly distributed random variable on $\Omega$ then the average of the sum of $f(U_i)\; i\in [1,\,N]$ converges to $I$ almost surely when $n$ tends to infinity, i.e.:
\begin{equation}
\int_{\Omega} f(x)\mathrm{d}x\approx \frac{1}{N}\sum_{i=1}^{N}f(U_i)
\end{equation}

Hence, while conventional numerical methods calculate the integrand at regularly spaced points, MC method samples the integrand at random points $U_i,\; i\in [1,N]$ ($N$ is the number of samples).\\

The critical issue with these points, is that they may not be equally distributed in the domain and this leads to the need to increase the number of samples, and, consequently, run-times.\\

This problem can be solved with QMC methods, making use of quasi-random numbers that are more well-distributed \cite{l2005recent}. Although quasi-random numbers come from a deterministic algorithm, they pass a statistical test of randomness.\\

Among  these methods those which make use of low discrepancy sequences (LDS)\cite{ortiz2008some} are based on the property of lack of an apparent pattern in the distance\footnote{6 is the most common separation distance up to about $n\approx 1.74\times 10^{35}$} between couples of primes and for this reason conforming a set of quasi-random numbers.\\

The application of the MC and QMC methods to the Hardy-Littlewood integral calculation has been
explored:

$$2C_2\int_{2}^{n}\frac{\mathrm{d}x}{(\ln\ln(x))^2}$$

using low discrepancy sequences (LDS) and Mathematica software \footnote{It is worth notice that a compensating constant $a \times 7.39$ has been used, depending on the limits of integration, the minimal and maximal values of the set of samples, and the dimensions of the integrand \cite{ortiz2008some}} (Annex I).\\
 The following Table (Fig.\ref{FigMCvsLDS}) provides the results of the MC and LDS methods:

%%%%
\begin{figure}[htb]
    \centering
    \includegraphics[width=1\textwidth]{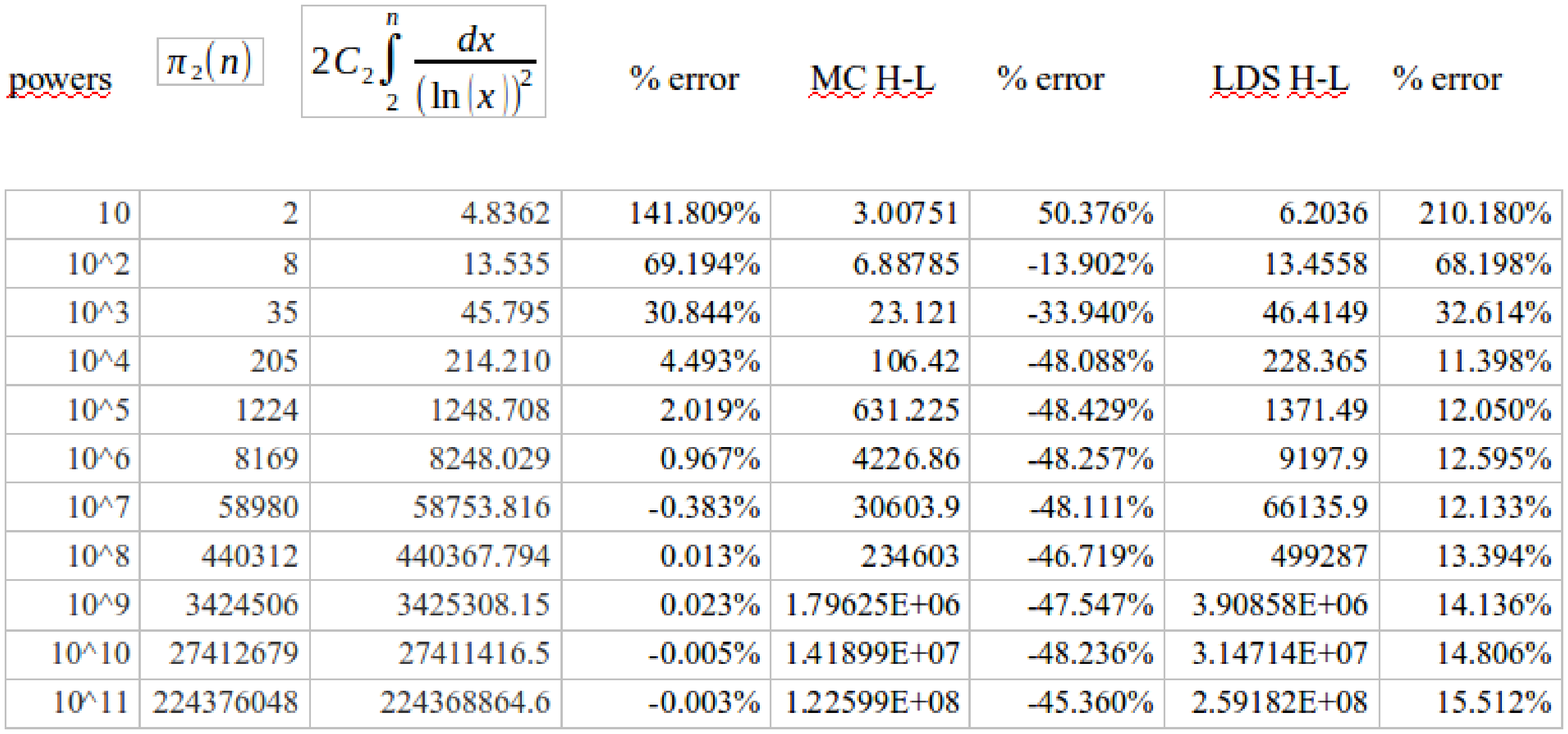}
\caption{Comparisons of MC and LDS methods for the first Hardy-Littlewood conjecture}\label{FigMCvsLDS}           
\end{figure}

Since the convergence rate of Monte Carlo method is close\footnote{It is rather slow: quadrupling the number of sampled points will halve the error} to $\mathcal{O\left(\frac{1}{\sqrt{N}}\right)}$, the error rate decreases as the value of $N$ increases (i.e. as a function $\pi_2(n)$ increases) as described in literature.\\

In the table shown in Fig.\ref{FigMCvsLDS}, the convergence  is not proved due to the low number of $N$ points  considered in the calculation ($n=10^{11},\; N=17548$), but the advantage of using LDS can be appreciated.\\

Finally, the following table (Fig. \ref{FigMC1011}) provides the results of the Monte Carlo method with a sufficient number of samples:

%%%%
\begin{figure}[htb]
    \centering
    \includegraphics[width=1\textwidth]{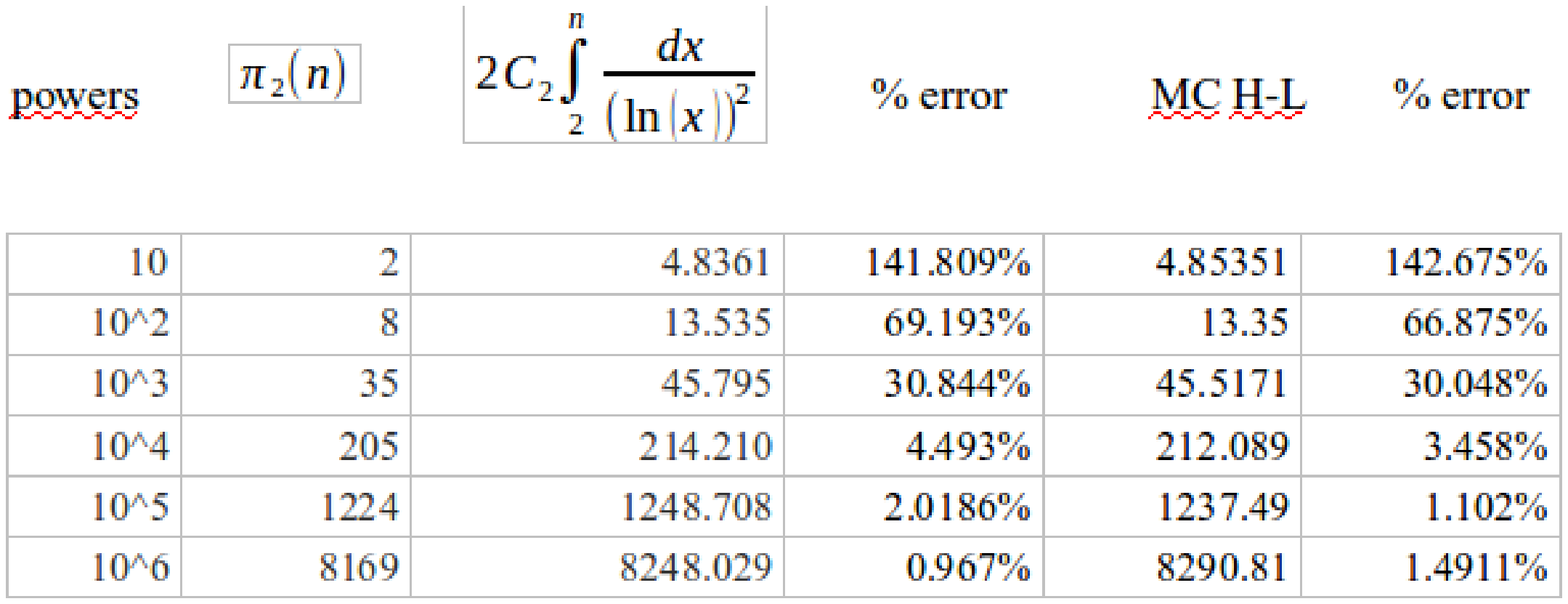}
\caption{Hardy-Littlewood conjecture using MC with a a larger set of $N$ points.}\label{FigMC1011}           
\end{figure}

\section{CONCLUSIONS}

In spite of many challenges and improvements due to numerical analysis, twin primes are still an unsolved problem in number theory. The first Hardy-Littlewood conjecture can be described as a milestone in this field.\\

This paper has proposed an empirical analysis of the twin primes distribution that leads to write the conjecture in terms of $\mathrm{mod}(10)$ classes marked by the same cardinality, according to a statistically random system.\\

Furthermore, starting from the conjecture, an elementary demonstration of the Brun’s theorem about the convergence of the sum of the reciprocal of the twin primes has been provided.\\

Finally, a less conventional method of calculation of the Hardy-Littlewood integral has been explored based on the MC and QMC methods involving the use of low discrepancy sequences (LDS).\\

 The result of the calculation with a sufficient number of samples is compelling and provides (for any given $n$ larger than $n=10^6$ say) a small relative error and an original example of application of these methods to the number theory.

\section{Annex:  MonteCarlo code  (Mathematica) }
{\color{teal}
\begin{verbatim}
powers = {10, 10^2, 10^3, 10^4, 10^5, 10^6, 10^7, 10^8, 10^9, 10^10, 
   10^11};
nooftwinp = {2, 8, 35, 205, 1224, 8169, 58980, 440312, 3424506, 
   27412679, 224376048};
HLconjecture = {4.8361883278, 13.5354875604, 45.7955004115, 
   214.2109398311, 1248.7087356371, 8248.0296898308, 58753.8164979342,
    440367.7942273770, 3425308.1557430851, 27411416.5322785837, 
   224368864.6811819439};

normalized = 2*ListMetadistances66/Max[ListMetadistances66];
(*This calls the list of LDS, named ListMetadistances66*)

mcHLintegrand = 
  Table[Table[
    1/(Log[x])^2, {x, 2, powers[[k]], (powers[[k]] - 2)/17547}], {k, 
    1, Length[powers]}];
 (*discretize integrand for calculation of integral for metadistances*)

mcHLsummatories = 
 Table[2 c2*
   Sum[powers[[k]]* mcHLintegrand[[k, i]]*
     RandomReal[]/Length[mcHLintegrand[[k]]], {i, 1, 
     Length[mcHLintegrand[[k]]] - 1}], {k, 1, Length[mcHLintegrand]}]
 (*calculate integral using MC, just up to the length of the discretized integrand*)

HLintegrand2 = 
  Table[1/(Log[x])^2, {x, 2, 
    powers[[2]], (powers[[2]] - 2)/(Length[ListMetadistances66])}];
  (*discretize the integrand for metadistances*)


summatories = 
 Table[2 c2*
   Sum[powers[[k]]*7.39* HLintegrand[[k, i]]*
     normalized[[i]]/Length[HLintegrand[[k]]], {i, 1, 
     Length[HLintegrand[[k]]] - 1}], {k, 1, Length[HLintegrand]}] 
     (*calculate all the integrals, for every upper limit of the integral (powers) *)

(*make the comparisons*)
comparisons = 
 Table[{pow2[[k]], pi2n[[k]], ScientificForm[HLconjecture[[k]], 3], 
   ScientificForm[(Abs[pi2n[[k]] - HLconjecture[[k]]])*100/pi2n[[k]], 
    3], ScientificForm[mcHLsummatories[[k]], 
    3], (Abs[pi2n[[k]] - mcHLsummatories[[k]]])*100/pi2n[[k]], 
   ScientificForm[summatories[[k]], 3], 
   N[(Abs[pi2n[[k]] - summatories[[k]]])*100/pi2n[[k]], 3]}, {k, 1, 
   Length[summatories]}]; PrependTo[comparisons, {"powers", 
  "\!\(\*SubscriptBox[\(\[Pi]\), \(2\)]\)(n)", "HL conj.", "% error", 
  "mc HL", "% error", "LDS HL", "% error"}]; MatrixForm[comparisons]     
\end{verbatim}
}

%%%%%%%%%%%%%%%%%%%%%%%%%%%

%\newpage

\bibliographystyle{siam} 
\bibliography{Twins.bib}

%\begin{thebibliography}{10}
%
%\bibitem{Eco}Eco, Umberto 2004. {\it History of Beauty}. Rizzoli, Italy.
%
%\bibitem{Fibonumbers} Lucas, Edouard. 1891. {\it Th\'eorie des nombres}. Paris: Gauthier-Villars et Fils, Impremeurs-Libraires.
%
%\bibitem{LucasNumbers} Niven, Ivan; Zuckerman, Herbert S.; Montgomery, Hugh L. 1991. {\it An Introduction to The theory of Numbers}. New York, John Wiley \& Sons, Inc. 
%
%\bibitem{GoldenBook}Hemenway,Priya. 2005. {\it Divine Proportion:  Phi in Art, Nature, and Science}. Sterling Publishing Co. pp. 127–129
%
%
%\end{thebibliography}

\end{document}